\documentclass[a4paper,12pt]{article}
\RequirePackage{amsthm,amsmath}

\usepackage{latexsym}
\usepackage{amsmath}
\usepackage{amsfonts}
\usepackage{amssymb}
\usepackage{psfrag}
\usepackage{graphicx}

\begin{document}

\theoremstyle{plain}
\newtheorem{thm}{Theorem}

\noindent\section*{The tempered discrete Linnik distribution} 
\noindent  
LUCIO BARABESI\\ {\it Department of Economics  \& Statistics, University of Siena, Italy}\ \
\vskip 0,2cm
\noindent
CAROLINA BECATTI\\ {\it Department of Economics  \& Statistics, University of Siena, Italy}\ \
\vskip 0,2cm
\noindent
MARZIA MARCHESELLI\\ {\it Department of Economics  \& Statistics, University of Siena, Italy}

\vskip0.6cm\noindent
{\bf ABSTRACT. A tempered version of the discrete Linnik distribution is introduced in order to obtain integer-valued distribution families connected to stable
laws. The proposal constitutes a generalization of the well-known
Poisson-Tweedie law, which is actually a tempered discrete stable law. The
features of the new tempered discrete Linnik distribution are explored by
providing a series of identities in law - which describe its genesis in terms
of mixture and compound Poisson law, as well as in terms of mixture discrete
stable law. A manageable expression of the corresponding probability function
is also provided and several special cases are analysed.}

\bigskip
\noindent {\it Key words}: L\'evy-Khintchine representation; Positive Stable distribution; Discrete Stable distribution; Mixture Poisson distribution; Compound Poisson distribution.\vskip0.5cm\noindent
{MSC 2010} Primary 62E10, secondary 62E15.
\vskip1cm
\noindent
{\bf 1. Introduction}
\vskip0.3cm
\noindent In recent years, heavy-tailed models - in primis, stable
distributions - have been used in a variety of fields, such as statistical
physics, mathematical finance and financial econometrics (see {\it e.g.} Rachev et
al., 2011, and references therein). However, these models may be partially
appropriate to provide a good fit to data, since their tails are too \lq\lq fat"  to describe empirical distributions, as remarked by Klebanov and\& Sl\'amov\'a (2015). In order to overcome this drawback, the so-called
\lq\lq tempered" versions of heavy-tailed distributions have been successfully
introduced (see {\it e.g.} Ros\'inski, 2007). Indeed, tempering allows for
models that are similar to original distributions in some central region,
even if possess lighter -  {\it i.e.} tempered - tails. Klebanov \& Sl\'amov\'a
(2015) have suitably discussed these issues and have suggested various
tempering techniques - emphasizing that tempering is not necessarily unambiguous.

In the framework of integer-valued laws, the Poisson-Tweedie is a well-known
tempered distribution introduced independently by Gerber (1991) and Hougaard
et al. (1997). This law is {\it de facto} the tempered counterpart of the Discrete
Stable law originally suggested by Steutel \& van Harn (1979). The
Poisson-Tweedie distribution encompasses classical families (such as the
Poisson), as well as large families (such as the Generalized Poisson Inverse
Gaussian and the Poisson-Pascal). Hence, this law may be very useful for
modelling data arising in a {\it plethora} of frameworks - for example, clinical
experiments (Hougaard, Lee \& Whitmore, 1997), environmental studies (El-Shaarawi, Zhu \& Joe, 2011) and scientometric analysis (Baccini, Barabesi \& Stracqualursi, 2016).

Christoph \& Schreiber (1998) emphasize that the Discrete Stable
distribution may be seen as the special case - for the limiting value of a
parameter - of the so-called Discrete Linnik distribution introduced by
Devroye (1993) and Pakes (1995). Therefore, the Discrete Linnik distribution
is more flexible than the Discrete Stable distribution and it could be very
useful to achieve its tempered version. In the present paper, we
preliminarily explore in detail the method - roughly outlined by Barabesi \&
Pratelli (2014a) - for obtaining integer-valued families of distributions
linked to stable and tempered stable laws. On the basis of these findings,
after revising some properties of the Discrete Linnik distribution, we
introduce its tempered counterpart. The new Tempered Discrete Linnik is
analysed thoroughly and its properties are given. In particular, some
stochastic representations for this law are obtained and the corresponding
probability function is achieved as a manageable finite sum.

The paper is organized as follows. In Sections 2 and 3, we revise and expand
the issues suggested by Barabesi \& Pratelli (2014a) for introducing
families of distributions connected to stable and tempered stable laws,
respectively. In Section 4, we survey the main features of the Discrete
Linnik distribution. Section 5 contains our proposal for the new Tempered
Discrete Linnik distribution, while in Section 6 we consider its main
properties. Finally, some conclusions are given in Section 7.

\vskip0.5cm
\noindent
{\bf 2. Integer-valued distribution families linked to stable laws}
\vskip0.15cm\noindent
Barabesi \& Pratelli (2014a) have suggested an approach - based on the definition of
subordinator - for devising integer-valued distribution families as mixtures
of stable laws, as well as tempered stable laws. In order to describe and to
develop at length their proposal, we first consider the absolutely-continuous
Positive Stable random variable (r.v.) - say $X_{PS}$ - with Laplace
transform given by
\begin{equation}
L_{X_{PS}}(t)=\exp(-\lambda t^\gamma)\text{ , Re}(t)>0,
\end{equation}
where $(\gamma,\lambda)\in]0,1]\times\mathbb{R}^{+}$ and
$\mathbb{R}^{+}=]0,\infty[$ (see {\it e.g.} Zolotarev, 1986, p.114). As is well
known, $\gamma$ is the so-called characteristic exponent -  {\it i.e.} a \lq\lq tail''
index - while $\lambda$ is actually a scale parameter. In order to emphasize
the dependence on $\gamma$ and $\lambda$, we eventually adopt the notation
$X_{PS}:=X_{PS}(\gamma,\lambda)$. 

We have also to introduce the integer-valued counterpart of the Positive
Stable r.v.,  {\it i.e.} the Discrete Stable r.v. $X_{DS}$ proposed by Steutel \&
van Harn (1979) with probability generating function (p.g.f.) given by
\begin{equation}
g_{X_{DS}}(s)=\exp(-\lambda(1-s)^\gamma)\text{ , }s\in[0,1],
\end{equation}
where in turn $(\gamma,\lambda)\in]0,1]\times\mathbb{R}^{+}$. For a survey of
the properties of this distribution, see {\it e.g.} Marcheselli, Barabesi \& Baccini (2008).
Similarly to the Positive Stable r.v., we also write
$X_{DS}:=X_{DS}(\gamma,\lambda)$. The parameter $\gamma$ is again a \lq\lq tail''
index, while $\lambda$ is a \lq\lq scale" parameter.

In order to clarify the meaning of the scale parameter for the Discrete
Stable law, we remind the \lq\lq thinning'' operator $\odot$ introduced by Steutel
\& van Harn (1979). If $X$ is an integer-valued r.v., the dot product
$\alpha\odot X$ is defined as
\begin{equation*}
\alpha\odot X=\sum_{i=1}^XZ_i\text{ ,}
\end{equation*}
where the $Z_i$'s are copies of Bernoulli r.v's of parameter $\alpha\in[0,1]$
independent of $X$. Obviously, the p.g.f. of $\alpha\odot X$ is given by $g_{\alpha\odot X}(s)=g_X(1-\alpha+\alpha s)$, where $g_X$ is the p.g.f. of
$X$. The dot product is also defined for $\alpha>1$, whenever $g_{\alpha\odot
X}$ is a proper p.g.f. (see Christoph \& Schreiber, 2001). Hence, $\lambda$
is a \lq\lq scale" parameter for the Discrete Stable law in the sense that
\begin{equation*}
X_{DS}(\gamma,\lambda)\overset{\mathcal{L}}{=}\lambda^{1/\gamma}\odot
X_{DS}(\gamma,1).
\end{equation*}
Let $\nu$ be a measure on $\mathbb{R}^{+}$ in such a way that
$\int\min(1,x)\nu($d$x)<\infty$. From the L\'evy-Khintchine representation
(see {\it e.g.} Sato, 1999, p.197) there exists a positive r.v. $Y$ with Laplace
transform given by
\begin{equation*}
L_Y(t)=\exp(-\eta\psi(t))\text{ , Re}(t)>0\text{ ,}
\end{equation*}
where $\eta\in\mathbb{R}^{+}$ and $\psi(t)=\int(1-\exp(-tx))\nu($d$x)$.
Moreover, let $X_P:=X_P(\lambda)$ represent a Poisson r.v. with parameter
$\lambda$,  {\it i.e.} the p.g.f. of $X_P(\lambda)$ is given by
$g_{X_P}(s)=\exp(-\lambda(1-s))$ with $s\in[0,1]$. Hence, if the r.v.'s $X_P$
and $Y$ are independent, the Mixture Poisson r.v.
\begin{equation*}
X_{MP}\overset{\mathcal{L}}{=}X_P(Y)
\end{equation*}
displays the p.g.f. given by
\begin{equation}
g_{X_{MP}}(s)=L_Y(1-s)=\exp(-\eta\psi(1-s))\text{ , }s\in[0,1].
\end{equation}
The Discrete Stable r.v. $X_{DS}(\gamma,\lambda)$ is obtained from expression
(3) when a stable subordinator is actually considered,  {\it i.e.}
$\nu($d$x)/$d$x\propto x^{-\gamma-1}I_{\mathbb{R}^{+}}(x)$ with
$\gamma\in]0,1[$ and where $I_B$ represents the indicator function of the set
$B$. Indeed, since the expression of $\nu$ gives rise to $\psi(t)\propto
t^\gamma$, the p.g.f. (2) promptly follows from expression (3) with
$\eta=\lambda$ and $\psi(t)=t^\gamma$. Moreover, since
$g_{X_{DS}}(s)=L_{X_{PS}}(1-s)$, it follows that
$Y\overset{\mathcal{L}}{=}X_{PS}(\gamma,\lambda)$ and
\begin{equation}
X_{DS}(\gamma,\lambda)\overset{\mathcal{L}}{=}X_P(X_{PS}(\gamma,\lambda)),
\end{equation}
which is actually equivalent to the identity in distribution emphasized by
Devroye (1993, Theorem in Section 1). On the basis of expression (4), a
general \lq\lq scale" mixture of Discrete Stable r.v.'s, say $X_{MDS}$, with a
mixturing absolutely-continuous positive r.v. $V$ having Laplace transform
$L_V$, may be achieved by considering the identity in distribution
\begin{equation}
X_{MDS}\overset{\mathcal{L}}{=}X_{DS}(\gamma,V)\overset{
\mathcal{L}}{=}X_P(X_{PS}(\gamma,V)),
\end{equation}
where the r.v.'s involved in the right-hand side are independent. Obviously,
(4) is achieved from (5) by assuming a degenerate distribution for $V$,  {\it i.e.} 
$P(V=\lambda)=1$. In addition, expression (5) may be also seen as the
stochastic \lq\lq scaling" of the Discrete stable law in terms of the operator
$\odot$,  {\it i.e.}
\begin{equation*}
X_{MDS}\overset{\mathcal{L}}{=}V^{1/\gamma}\odot
X_{DS}(\gamma,1)\overset{\mathcal{L}}{=}V^{1/\gamma}\odot
X_P(X_{PS}(\gamma,1))\,.
\end{equation*}
Moreover, from (5), it is apparent that the p.g.f. of the r.v. $X_{MDS}$
turns out to be
\begin{equation}
g_{X_{MDS}}(s)=L_V((1-s)^\gamma)\text{ , }s\in[0,1].
\end{equation}
Hence, families of mixture of Discrete Stable r.v.'s can be generated on the
basis of (5) and (6) by suitably selecting the r.v. $V$.

We conclude with a final remark on the p.g.f. (6). Let $X_S$ be a Sibuya r.v.
(as named by Devroye, 1993) with p.g.f.
\begin{equation*}
g_{X_S}(s)=1-(1-s)^\gamma\text{ , }s\in[0,1],
\end{equation*}
where $\gamma\in]0,1]$ (for a recent survey of this law, see Huillet, 2016).
Indeed, the Sibuya distribution is a special case of the (shifted) Negative
Binomial Beta distribution introduced by Sibuya (1979) with parameters given
by $1$, $\gamma$ and $(1-\gamma)$. In the following, the Sibuya r.v. is also
denoted by $X_S:=X_S(\gamma)$. Therefore, expression (6) may be also
interestingly reformulated as
\begin{equation}
g_{X_{MDS}}(s)=L_V(1-g_{X_S}(s))\text{ , }s\in[0,1].
\end{equation}
Thus, if $L_V$ displays a suitable structure, expression (7) eventually gives
rise to a representation of the r.v. $X_{MDS}$ in terms of a compound r.v.
with a compounding Sibuya r.v. As a quite easy example, the r.v. $X_{DS}$ may
be also expressed as a compound Poisson r.v. as 
\begin{equation*}
g_{X_{DS}}(s)=\exp(-\lambda(1-g_{X_S}(s)))\text{ , }s\in[0,1],
\end{equation*}
and hence
\begin{equation*}
X_{DS}(\gamma,\lambda)\overset{\mathcal{L}}{=}\sum_{i=1}^ZW_i,
\end{equation*}
where $Z\overset{\mathcal{L}}{=}X_P(\lambda)$ and the $W_i$'s are i.i.d. r.v.'s such that $W_i\overset{\mathcal{L}}{=}X_S(\gamma)$ - which are in turn
independent of $Z$.
\vskip0.5cm
\noindent {\bf 3. Integer-valued distribution families linked to tempered stable laws}
\vskip0.15cm
\noindent
First, for subsequent use, we provide some issues on the so-called Tweedie
distribution (for a recent survey of this law, see Barabesi et al., 2016).
The Tweedie distribution is actually a Tempered Positive Stable distribution
introduced by Hougaard (1986). Hence, for this reason, in the following we
denote the Tweedie r.v. as $X_{TPS}$. With a slight change in the
parameterization proposed by Hougaard (1986), the Laplace transform of the
r.v. $X_{TPS}$ is given by
\begin{equation}
L_{X_{TPS}}(t)=\exp(\text{sgn}(\gamma)\lambda(\theta^\gamma-(\theta+t)^\gamma))\text{ , Re}(t)>0,
\end{equation}
where
$(\gamma,\lambda,\theta)\in\{]-\infty,1]\times\mathbb{R}^{+}\times
\mathbb{R}^{+}\}\cup\{]0,1]\times\mathbb{R}^{+}\times\{0\}\}$. The
formulation proposed in expression (8) is convenient, since it avoids to
define the Laplace transform for analytical continuity for $\gamma=0$ as in
the case of the parameterization considered by Hougaard (1986). Moreover, it
is worth noting that $\theta$ is actually the \lq\lq tempering" parameter. This
is at once apparent for $\gamma\in]0,1]$ from the following identity
\begin{equation*}
L_{X_{TPS}}(t)=\frac{L_{X_{PS}}(\theta+t)}{L_{X_{PS}}(\theta)},
\end{equation*}
which reveals the exponential \lq\lq nature" of the tempering. By following the
usual route, we also write $X_{TPS}:=X_{TPS}(\gamma,\lambda,\theta)$.
Obviously,  for
$\gamma\in]0,1]$ it holds
$$X_{TPS}(\gamma,\lambda,0)\overset{\mathcal{L}}{=}X_{PS}(\gamma,\lambda).$$
It should be strongly remarked that the tempering extends the range of
parameter values (with respect to the Positive Stable distribution) for the
parameter $\gamma$ - which may assume negative values, even if $\theta$ must
be strictly positive in such a case. This is an interesting feature, since
for $\gamma\in\mathbb{R}^{-}$ where $\mathbb{R}^{-}=]-\infty,0[$, it is
immediate to reformulate the r.v. $X_{TPS}$ as a compound Poisson of Gamma
r.v.'s (see {\it e.g.} Barabesi et al., 2016). More precisely, let the r.v.
$X_G:=X_G(\lambda,\delta)$ be distributed according to a Gamma law with
corresponding Laplace transform given by $L_{X_G}(t)=(1+\lambda t)^{-\delta}$
with Re$(t)>0$ and $(\lambda,\delta)\in\mathbb{R}^{+}\times\mathbb{R}^{+}$
(obviously, $\lambda$ is the scale parameter and $\delta$ is the shape
parameter). Hence, on the basis of (8) and owing to the reproductive property
of the Gamma distribution with respect to the shape parameter, the following
identity in distribution holds for $\gamma\in\mathbb{R}^{-}$
\begin{equation}
X_{TPS}(\gamma,\lambda,\theta)\overset{\mathcal{L}}{=}X_G(1/\theta,-\gamma
X_P(\lambda\theta^\gamma)).
\end{equation}
Hence, in this case the r.v. $X_{TPS}$ displays a mixed distribution, given
by a convex combination of a Dirac distribution (with mass at zero) and an
absolutely-continuous distribution (a very useful property for modelling data
with an excess of zeroes, see Barabesi et al., 2016).

In order to extend the issues of Section 2, we introduce a tempered stable subordinator,  {\it i.e.} $$\nu(dx)/dx\propto\exp(-\theta
x)x^{-\gamma-1}I_{\mathbb{R}^{+}}(x),$$ where
$(\gamma,\theta)\in\{]-\infty,1[
\times\mathbb{R}^{+}\}\cup\{]0,1[\times\{0\}\}$. 
In such a case, on the basis of (3) the so-called Poisson-Tweedie
distribution - which is actually a Tempered Discrete Stable distribution - is
achieved (for more about the Poisson-Tweedie law, see Baccini, Barabesi \& Stracqualursi, 2016,
and El-Shaarawi, Zhu \& Joe, 2011). Indeed, if $X_{TDS}$ denotes a Tempered
Discrete Stable r.v., since the expression of $\nu$ provides
$\psi(t)\propto(\theta+t)^\gamma-\theta^\gamma$, from (3) with $\eta=\lambda$
and $\psi(t)=(\theta+t)^\gamma-\theta^\gamma$ the following p.g.f. is obtained
\begin{equation*}
g_{X_{TDS}}(s)=\exp(\text{sgn}(\gamma)\lambda(\theta^\gamma-(\theta+1-s)^%
\gamma))\text{ , }s\in[0,1].
\end{equation*}
Hence - as expected - it holds $g_{X_{TDS}}(s)=L_{X_{TPS}}(1-s)$. However, in
order to be consistent with the existing literature and for practical
convenience, we prefer to reparameterize the previous p.g.f. by assuming that
$\gamma=a$, $\lambda=bc^a$ and $\theta=1/c-1$, in such a way that 
\begin{equation}
g_{X_{TDS}}(s)=\exp(\text{sgn}(a)b((1-c)^a-(1-cs)^a)), s\in[0,1],
\end{equation}
where
$(a,b,c)\in\{]-\infty,0]\times\mathbb{R}^{+}\times[0,1[\}\cup\{]0,1]\times
\mathbb{R}^{+}\times[0,1]\}$. The p.g.f. (10) is provided as a slight
modification of the formulation suggested by El-Shaarawi et al. (2011).
Indeed, the parameterization considered in (10) avoids to define the p.g.f.
for analytical continuity for $a=0$. In this special case, $X_{TDS}$
degenerates at $0$,  {\it i.e.} $P(X_{TDS}=0)=1$. As usual, we write
$X_{TDS}:=X_{TDS}(a,b,c)$.

It is worth noting that $c$ actually represents the \lq\lq tempering" parameter.
Indeed, for $a\in]0,1]$ and by considering the r.v. $X_{DS}(a,b)$, the
following identity
\begin{equation*}
L_{X_{TDS}}(t)=\frac{g_{X_{DS}}(cs)}{g_{X_{DS}}(c)}
\end{equation*}
emphasizes the geometric \lq\lq nature" of the tempering. In addition, it
straightforwardly holds $X_{DS}(a,b)\overset{\mathcal{L}}{=}X_{TDS}(a,b,1)$
for $a\in]0,1]$. In turn, similarly to the Tweedie distribution, tempering
extends the range of parameter values (with respect to the Discrete Stable
distribution) for the parameter $a$ - which may assume negative values, even
if $c$ must be strictly less than unity in such a case. Finally, from (3)
with $Y\overset{\mathcal{L}}{=}X_{TPS}(a,bc^a,1/c-1)$, it holds 
\begin{equation}
X_{TDS}(a,b,c)\overset{\mathcal{L}}{=}X_P(X_{TPS}(a,bc^a,1/c-1)),
\end{equation}
which actually constitutes the identity in distribution remarked by Hougaard
et al. (1997) and which generalizes expression (4).

On the basis of expression (11), it is at once apparent that a general
\lq\lq scale" mixture of Tempered Discrete Stable r.v.'s, say $X_{MTDS}$, with a
mixturing absolutely-continuous positive r.v. $V$ having Laplace transform
$L_V$, may be achieved by considering the identity in distribution - which
generalizes (5) - given by
\begin{equation}
X_{MTDS}\overset{\mathcal{L}}{=}X_{TDS}(a,V,c)\overset{%
\mathcal{L}}{=}X_P(X_{TPS}(a,c^aV,1/c-1)),
\end{equation}
where the r.v.'s involved in the right-hand side are independent. Obviously,
(11) is achieved from (12) by assuming a degenerate distribution for $V$,
 {\it i.e.} $P(V=b)=1$. Moreover, the corresponding p.g.f. turns out to be
\begin{equation}
g_{X_{MTDS}}(s)=L_V(\text{sgn}(a)((1-cs)^a-(1-c)^a))\text{ , }s\in[0,1].
\end{equation}
Hence, families of mixture of Tempered Discrete Stable r.v.'s can be
generated by means of (12) - and accordingly (13) - by suitably selecting the
r.v.$\,V$. In order to reformulate (13) similarly to (7) when $a\in]0,1]$, let
$X_{GDS}:=X_{GDS}(\gamma,\tau)$ be the Geometric Down-weighting Sibuya r.v.
(introduced by Zhu \& Joe, 2009) with the p.g.f.
\begin{equation*}
g_{X_{GDS}}(s)=1-g_{X_S}(\tau)+g_{X_S}(\tau s)=1+(1-\tau)^\gamma-(1-\tau
s)^\gamma,
\end{equation*}
where $(\gamma,\tau)\in]0,1]\times]0,1]$. Hence, by considering the Geometric
Down-weighting Sibuya r.v. $X_{GDS}(a,c)$, for $a\in]0,1]$ expression (13)
may be also interestingly rewritten as
\begin{equation}
g_{X_{MTDS}}(s)=L_V(1-g_{X_{GDS}}(s))\text{ , }s\in[0,1].
\end{equation}
Moreover, let $X_{NB}:=X_{NB}(\pi,\delta)$ be a Negative Binomial r.v. with
p.g.f. given by 
\begin{equation*}
g_{X_{NB}}(s)=\left(\frac{1-\pi}{1-\pi
s}\right)^\delta=(1+\lambda(1-s))^{-\delta}, s\in[0,1],
\end{equation*}
where $(\pi,\delta)\in]0,1[\times\mathbb{R}^{+}$, while
$\lambda=\pi/(1-\pi)\in\mathbb{R}^{+}$. Thus, when $a\in\mathbb{R}^{-}$ and
by considering the Negative Binomial r.v. $X_{NB}(c,-a)$, expression (13) may
be also rewritten as
\begin{equation}
g_{X_{MTDS}}(s)=L_V((1-c)^a(1-g_{X_{NB}}(s)))\text{ , }s\in[0,1].
\end{equation}
Therefore, if $L_V$ displays a suitable structure, expression (14) and (15)
eventually gives rise to a representation of the r.v. $X_{MDS}$ in terms of a
compound r.v. with compounding Sibuya and Negative Binomial r.v.'s,
respectively.
\vskip0.5cm\noindent
{\bf 4. Some remarks on the genesis and properties of the discrete Linnik distribution}
\vskip0.15cm\noindent By following the formulation adopted by Christoph \& Schreiber
(1998), the p.g.f. of the r.v. $X_{DL}$ distributed according to the Discrete
Linnik law is given by
\begin{equation}
g_{X_{DL}}(s)=(1+\lambda(1-s)^\gamma/\delta)^{-\delta}\text{ ,
}s\in[0,1],
\end{equation}
where
$(\gamma,\lambda,\delta)\in]0,1]\times\mathbb{R}^{+}\times\mathbb{R}^{+}$.
For a detailed description of the main features of the law, see Christoph \&
Schreiber (1998, 2001). It is apparent that the p.g.f. of the Discrete Stable
r.v. $X_{DS}$ is achieved as $\delta\rightarrow\infty$. In addition, it
should be remarked that the Discrete Linnik law is defined for some negative
$\delta$ in such a way that $\lambda\leq|\delta|(1-\gamma)$ and, in this
case, the distribution is also named Generalized Sibuya (Huillet, 2016).
However, the results of this Section are mainly given for positive $\delta$.
Thus, it is apparent that the Discrete Linnik law is very flexible and may
encompass a large variety of distribution families ranging from light-tailed
laws ({\it e.g.} the Binomial law for suitable $\lambda$ and negative integer
$\delta$ when $\gamma=1$) to heavy-tailed laws ({\it e.g.} the Discrete
Mittag-Leffler law for $\delta=1$). 

In order to emphasize the dependence on the parameters, we also adopt the
notation $X_{DL}:=X_{DL}(\gamma,\lambda,\delta)$. It should be remarked that
the parameter $\gamma$ is in turn a \lq\lq tail'' index and $\lambda$ is a
\lq\lq scale" parameter, while $\delta$ is a shape parameter. Indeed, $\lambda$
is a \lq\lq scale" parameter for the Discrete Linnik law in the sense that 
\begin{equation*}
X_{DL}(\gamma,\lambda,\delta)\overset{\mathcal{L}}{=}\lambda^{1/\gamma}\odot
X_{DL}(\gamma,1,\delta).
\end{equation*}
For $\delta\in\mathbb{R}^{+}$ and on the basis of the remarks provided in
Section 2, it is promptly proven that the Discrete Linnik r.v. is a \lq\lq scale"
mixture of Discrete Stable r.v.'s obtained by selecting
$V\overset{\mathcal{L}}{=}X_G(\lambda/\delta,\delta)$ in (5). In this case,
the p.g.f. (16) is obtained by means of (6), while from (5) it also holds
\begin{equation}
X_{DL}(\gamma,\lambda,\delta)\overset{\mathcal{L}}{=}X_{DS}(\gamma,X_G(%
\lambda/\delta,\delta))\overset{\mathcal{L}}{=}X_P(X_{PS}(\gamma,X_G(\lambda/%
\delta,\delta))),
\end{equation}
which is actually similar to the identity in distribution obtained by Devroye
(1993, Section 2). A further remark is in order, since expression (17) may be
also rephrased in terms of the absolutely-continuous Positive Linnik r.v.
$X_{PL}:=X_{PL}(\gamma,\lambda,\delta)$ with Laplace transform given by
\begin{equation*}
L_{X_{PL}}(t)=(1+\lambda t^\gamma/\delta)^{-\delta}\text{ , Re}(t)>0,
\end{equation*}
where
$(\gamma,\lambda,\delta)\in]0,1]\times\mathbb{R}^{+}\times\mathbb{R}^{+}$ -
according to the formulation provided by Christoph \& Schreiber (2001). The
absolutely-continuous Positive Linnik has been introduced by Pakes (1995)
after that Linnik (1962, p.67) had identified its symmetric version in a
simple case. Obviously, the Laplace transform (1) of the Positive Stable r.v.
$X_{PS}(\gamma,\lambda)$ is achieved as $\delta\rightarrow\infty$. For more
details on the Positive Linnik distribution, see {\it e.g.} Jose et al. (2010).
Since it easily proven that
\begin{equation*}
X_{PL}(\gamma,\lambda,\delta)\overset{\mathcal{L}}{=}X_{PS}(\gamma,X_G(%
\lambda/\delta,\delta)),
\end{equation*}
where the r.v.'s involved in the right-hand side are independent, by means of
expression (17) it promptly follows that
\begin{equation}
X_{DL}(\gamma,\lambda,\delta)\overset{\mathcal{L}}{=}X_P(X_{PL}(\gamma,%
\lambda,\delta)),
\end{equation}
which actually generalizes (4). Indeed, expression (4) is recovered from
expression (18) as $\delta\rightarrow\infty$. Finally, we remark that the
identity in distribution (17) is very suitable for random variate generation.
Indeed, many generators for Poisson and Gamma variates are available in
statistical literature, while Positive Stable variates are readily obtained
by means of the well-known Kanter's representation (Kanter, 1975). A similar
remark applies to identity (4).

For $\delta\in\mathbb{R}^{+}$, the r.v. $X_{DL}$ may be also expressed as a
compound Negative Binomial r.v. with a compounding Sibuya r.v. Indeed, from
expression (7) with $V\overset{\mathcal{L}}{=}X_G(\lambda/\delta,\delta)$, it
follows that
\begin{equation*}
g_{X_{DL}}(s)=(1+\lambda(1-g_{X_S}(s))/\delta)^{-\delta}\text{ ,
}s\in[0,1]\text{ ,}
\end{equation*}
and hence
\begin{equation*}
X_{DL}(\gamma,\lambda,\delta)\overset{\mathcal{L}}{=}\sum_{i=1}^ZW_i,
\end{equation*}
where $Z\overset{\mathcal{L}}{=}X_{NB}(\lambda/(\delta+\lambda),\delta)$ and
the $W_i$'s are i.i.d. r.v.'s such that
$W_i\overset{\mathcal{L}}{=}X_S(\gamma)$ - which are in turn independent of
$Z$.
\vskip0.5cm\noindent {\bf 5. The tempered discrete Linnik distribution}
\vskip0.15cm\noindent
On the basis of the issues discussed in Sections 2, 3 and 4, we are ready to introduce a tempered
version of the Discrete Linnik distribution. Indeed, in a complete
parallelism with Section 4, the Tempered Discrete Linnik r.v. may be given as
a \lq\lq scale" mixture of the Tempered Discrete Stable r.v.'s by selecting
$V\overset{\mathcal{L}}{=}X_G(bd,1/d)$ in (12). In this case, if the Tempered
Discrete Linnik r.v. is denoted by $X_{TDL}$, from expression (13) the
corresponding p.g.f. turns out to be
\begin{equation}
g_{X_{TDL}}(s)=(1+\text{sgn}(a)bd((1-cs)^a-(1-c)^a))^{-1/d}\text{ ,
}s\in[0,1],
\end{equation}
with $(a,b,c,d)\kern-1mm\in\kern-1mm\{]-\infty,0]\times\mathbb{R}^{+}\times[0,1[\times\mathbb{R}^{+}
\}\cup\{]0,1]\times\mathbb{R}^{+}\times[0,1]\times\mathbb{R}^{+}\}.$ Several
comments on the parameterization are in order. First, similarly to the
Discrete Linnik law, the Tempered Discrete Linnik law is also defined for
some negative $\delta$ - a distribution that could be named Tempered
Generalized Sibuya by extending the definition by Huillet (2016) - even if we
mainly deal with positive $d$. Secondly, the adopted formulation allows to
achieve the Poisson-Tweedie distribution as $d$ approaches zero - and this
could be preferable with respect to a limit at infinity. Moreover, for $a=0$
the r.v. $X_{TDL}$ degenerates at $0$. Thus, the p.g.f. (19) seems a natural
generalization of the p.d.f. (10). As usual, we also adopt the notation
$X_{TDL}:=X_{TDL}(a,b,c,d)$. Obviously, it promptly holds
$X_{DL}(a,b,d)\overset{\mathcal{L}}{=}X_{TDL}(a,b,1,d)$ for $a\in]0,1]$.
Analogously to the Poisson-Tweedie distribution, it is worth noting that
tempering extends the range of parameter values for the parameter $a$ with
respect to the Discrete Linnik distribution. 

For $d\in\mathbb{R}^{+}$, similarly to (11) and by means of (12), the
following identity in distribution holds
\begin{equation}
\begin{aligned}[t]
X_{TDL}(a,b,c,d)&\overset{\mathcal{L}}{=}X_{TDS}(a,X_G(bd,1/d),c))\\
&\overset{\mathcal{L}}{=}X_P(X_{TPS}(a,X_G(bdc^a,1/d),1/c-1)
\end{aligned}
\end{equation}
which actually generalizes (17). We remark that the identity (20) is suitable
for random variate generation, if Tempered Positive Stable variates are
available. Such a generation is straightforward on the basis of (9) if
$a\in\mathbb{R}^{-}$, since in this case
\begin{equation}
X_{TDL}(a,b,c,d)\overset{\mathcal{L}}{=}X_P(X_G(c/(1-c),-aX_P(X_G(bd(1-c)^a,
1/d)))).
\end{equation}
In contrast, for $a\in]0,1]$ the algorithms proposed by Barabesi et al. (2016) and Devroye (2009) for the generation of Tempered Positive Stable
variates may be considered. In any case, Barabesi \& Pratelli (2014b, 2015)
have suggested algorithms based on a rejection method linked to the p.g.f.,
which may be suitable for the direct generation of Tempered Discrete Linnik
variates. Actually, Baccini, Barabesi \& Stracqualursi (2016) have introduced a similar proposal
for the Tempered Discrete Stable law.

It should be remarked that - on the basis of (3) - an absolutely-continuous
Tempered Positive Linnik r.v.
$X_{TPL}:=X_{TPL}(\gamma,\lambda,\theta,\delta)$ may be given with Laplace
transform given by
\begin{equation*}
L_{X_{TPL}}(t)=(1+\text{sgn}(\gamma)\lambda((\theta+t)^\gamma-\theta^\gamma)/
\delta)^{-\delta}\text{ , Re}(t)>0,
\end{equation*}
where
$(\gamma,\lambda,\theta,\delta)\kern-1mm\in\kern-1mm\{]-\infty,1]\times\mathbb{R}^{+}\times
\mathbb{R}^{+}\times\mathbb{R}^{+}\}\cup\{]0,1]\times\mathbb{R}^{+}\times\{0
\}\times\mathbb{R}^{+}\}$. In turn, the Laplace transform (8) of the Tempered
Positive Stable r.v. $X_{TPS}(\gamma,\lambda)$ is achieved as
$\delta\rightarrow\infty$. Since it holds that
\begin{equation*}
X_{TPL}(\gamma,\lambda,\theta,\delta)\overset{\mathcal{L}}{=}X_{TPS}(\gamma,
X_G(\lambda/\delta,\delta),\theta),
\end{equation*}
where the r.v.'s involved in the right-hand side are independent, from
expression (20) it promptly follows that
\begin{equation*}
X_{TDL}(a,b,c,d)\overset{\mathcal{L}}{=}X_P(X_{TPL}(a,b,c,1/d)),
\end{equation*}
which actually generalizes (18).

For $d\in\mathbb{R}^{+}$, on the basis of (14) and (15), the r.v. $X_{TDL}$
may be also expressed as a compound distribution. Indeed, when $a\in]0,1]$
from (14) it holds
\begin{equation*}
g_{X_{TDL}}(s)=(1+bd(1-g_{X_{GDS}}(s)))^{-1/d}, s\in[0,1],
\end{equation*}
and hence
\begin{equation}
X_{TDL}(a,b,c,d)\overset{\mathcal{L}}{=}\sum_{i=1}^ZW_i,
\end{equation}
where $Z\overset{\mathcal{L}}{=}X_{NB}(bd/(1+bd),1/d)$ and the $W_i$'s are
i.i.d. r.v.'s such that $W_i\overset{\mathcal{L}}{=}X_{GDS}(a,c)$ - which are
in turn independent of $Z$. Moreover, when $a\in\mathbb{R}^{-}$ expression
(15) may be rewritten as
\begin{equation*}
g_{X_{TDL}}(s)=(1+bd(1-c)^a(1-g_{X_{NB}}(s)))^{-1/d}, s\in[0,1],
\end{equation*}
and hence (22) holds with
$Z\overset{\mathcal{L}}{=}X_{NB}(bd(1-c)^a/(1+bd(1-c)^a),1/d)$ and the
$W_i$'s are i.i.d. r.v.'s such that $W_i\overset{\mathcal{L}}{=}X_{NB}(c,-a)$
- which are in turn independent of $Z$. Owing to the reproductive property of
the Negative Binomial law, expression (22) also provides 
\begin{equation}
X_{TDL}(a,b,c,d)\overset{\mathcal{L}}{=}X_{NB}(c,
-aX_{NB}(bd(1-c)^a/(1+bd(1-c)^a),1/d)).
\end{equation}
It is worth remarking that identities (21) and (23) are straightforwardly
equivalent since
$X_{NB}(\pi,\delta)\overset{\mathcal{L}}{=}X_P(X_G(\pi/(1-\pi),\delta))$. 
\vskip0.5cm
\noindent{\bf 6.Further features of the tempered discrete Linnik distribution}
\vskip0,15cm\noindent
First, in order to obtain a closed form for the probability function (p.f.) of the
Tempered Discrete Linnik r.v., we provide a general result on the p.f. of
r.v.'s displaying a general type of p.g.f. - which encompasses (13) as
special case. Indeed, let us suppose that the r.v. $X$ has p.g.f. given by
\begin{equation*}
g_X(s)=\varphi(\alpha+\beta(1-\phi s)^\gamma)\text{ , }s\in[0,1],
\end{equation*}
where $\gamma$, $\alpha$, $\beta$ and $\phi$ are real constants, in such a
way that $\phi\in[0,1]$. If $\varphi$ is analytic in a neighbourhood of
$(\alpha+\beta)$, it follows that
\begin{equation*}
\begin{aligned}[t]
P(X=k)&=\frac{1}{k!}\frac{d^k}{ds^k}\left.\varphi(\alpha+\beta+\beta((1-\phi
s)^\gamma-1)\right|_{s=0}\\
&=\frac{1}{k!}\,\sum_{m=0}^k\frac{\beta^m}{m!}\frac{d^m}{ds^m}\left.
\varphi(s)\right|_{s=\alpha+\beta}\,\frac{d^k}{ds^k}\left.((1-\phi
s)^\gamma-1)^m\right|_{s=0}\\
&=\frac{1}{k!}\,\sum_{m=0}^k\frac{\beta^m}{m!}\frac{d^m}{ds^m}\left.
\varphi(s)\right|_{s=\alpha+\beta}\,\sum_{j=0}^m(-1)^{m-j}\binom{m}{j}
\frac{d^k}{ds^k}\left.(1-\phi s)^{j\gamma}\right|_{s=0}\\
&=(-\phi)^k\,\sum_{m=0}^k\frac{(-\beta)^m}{m!}\frac{d^m}{ds^m}\left.
\varphi(s)\right|_{s=\alpha+\beta}\,C_{\gamma,m}(k),
\end{aligned}
\end{equation*}
where
\begin{equation*}
C_{\gamma,m}(k)=\sum_{j=0}^m(-1)^j\binom{m}{j}\binom{\gamma j}{k}.
\end{equation*}
Thus, the p.d.f. of $X$ may be achieved as a finite sum. This result is very
useful for obtaining the p.f. of the r.v. $X_{TDL}$ by setting
$\varphi(x)=x^{-1/d}$ and the values $\gamma=a$,
$\alpha=1-$sgn$(a)bd(1-c)^a$,$\,\beta=$sgn$(a)bd$ and $\phi=c$,  {\it i.e.}
\begin{equation}
P(X_{TDL}=k)=(-c)^k\sum_{m=0}^k\binom{-1/d}{m}\frac{(-\text{sgn}(a)bd)^m}{(1+
\text{sgn}(a)bd(1-(1-c)^a))^{m+1/d}}\,C_{a,m}(k)\text{ \ \ }
\end{equation}
for $k\in\mathbb{N}$. Similarly, the p.f. of Tempered Discrete Stable r.v. is
achieved by setting $\varphi(x)=\exp(x)$ and the values $\gamma=a$,
$\alpha=$sgn$(a)b(1-c)^a$,$\,\beta=-$sgn$(a)b$ and $\phi=c$,  {\it i.e.}
\begin{equation}
P(X_{TDS}=k)=\exp(-\text{sgn}(a)b(1-(1-c)^a))(-c)^k\sum_{m=0}^k\frac{(
\text{sgn}(a)b)^m}{m!}\,C_{a,m}(k)
\end{equation}
for $k\in\mathbb{N}$. It is worth remarking that (25) is a
ctually equivalent
to the expression provided by Baccini, Barabesi \& Stracqualursi (2016).

On the basis of these findings, a series of comments on some particular
values of the parameters of the Tempered Discrete Linnik distribution for
$d\in\mathbb{R}^{+}$ are worthwhile. Indeed, when $a=1$, expression (19)
provides $X_{TDL}(1,b,c,d)\overset{\mathcal{L}}{=}X_{NB}(bcd/(1+bcd),1/d)$
and hence the r.v. $X_P(bc)$ is obtained as $d\rightarrow 0$. It should be
noticed that the Negative Binomial is also achieved when $a\in]0,1]$ and
$bd(1-c)^a=1$. Moreover, for $a=0$, the r.v. $X_{TDL}(0,b,c,d)$ has a
degenerate distribution at zero, as previously remarked. In addition, for
$d=1$ a tempered version of the Discrete Mittag-Leffler is achieved.

For $a=1/2$, the Tempered Discrete Linnik law is actually a generalization of
the Poisson Inverse Gaussian law (for a discussion of this law, see {\it e.g.}
Johnson et al., 2005, p.484), which is obtained as $d\rightarrow 0$. Thus,
this new distribution could be suitably named Negative Binomial Inverse
Gaussian. Moreover, in this special case expression (24) remarkably reduce to
a single sum, since
\begin{equation*}
C_{1/2,m}(k)=(-1)^k\,\frac{2^{-2k+m}m}{2k-m}\binom{2k-m}{k}
\end{equation*}
for $k\in\mathbb{N}^{+}$, while $C_{1/2,m}(0)=1$. In addition, since it holds
\begin{equation*}
C_{1/2,m}(k+1)=-\frac{2k-m}{4(2k-m+2)}\,C_{1/2,m}(k),
\end{equation*}
the probabilities in (24) may be sequentially evaluated - so that the actual
computational burden is small. 

When $a=-1$, the Tempered Discrete Linnik is in turn a generalization of the
Poly\'a-Aeppli law (for more about this distribution, see {\it e.g.} Johnson et
al., 2005, p.410), which is obtained as $d\rightarrow 0$. In this case also,
expression (24) and (25) remarkably reduce to a single sum, since
\begin{equation*}
C_{-1,m}(k)=(-1)^{k+m}\binom{k-1}{m-1}
\end{equation*}
for $k\in\mathbb{N}^{+}$, while $C_{-1,m}(0)=1$. Similarly to the case
$a=1/2$, the probabilities in (24) may be sequentially evaluated, since it
holds
\begin{equation*}
C_{-1,m}(k+1)=-\frac{m}{k}\,C_{-1,m}(k).
\end{equation*}

As to the main descriptive indexes, on the basis of (19) and after tedious
algebra, it follows that the expectation and the variance of the r.v.
$X_{TDL}$ are respectively given by
\begin{equation*}
\mu={\rm E}[X_{TDL}]=\text{sgn}(a)abc(1-c)^{a-1}
\end{equation*}
and
\begin{equation*}
\sigma^2={\rm Var}[X_{TDL}]=d\mu^2+\frac{(1-ac)\mu}{1-c}.
\end{equation*}
It is worth noting that $\mu$ does not depend on the parameter $d$ -  {\it i.e.},
the Tempered Discrete Linnik r.v. and the Poisson-Tweedie r.v. actually
display the same expectation. However, since the dispersion index is given by
\begin{equation*}
D=\frac{\sigma^2}{\mu}=d\mu+\frac{1-ac}{1-c},
\end{equation*}
the Poisson-Tweedie distribution may solely display over-dispersion (since
$D\geq 1$ when $d=0$), while the Tempered Discrete Linnik distribution may
accommodate for under-dispersion (for some admissible values $d\,<\,0$), as
well as for over-dispersion (when $d>0$). Hence, the Tempered Discrete Linnik
law substantially extend the range of the dispersion index with respect to
the Poisson-Tweedie law.

After further tedious algebra, it also follows
\begin{equation*}
m_3={\rm E}[(X_{TDL}-\mu)^3]=\frac{\sigma^4}{\mu}+d\mu\sigma^2+
\frac{c(1-a)\mu}{(1-c)^2}
\end{equation*}
and
\begin{equation*}
m_4={\rm E}[(X_{TDL}-\mu)^4]=3(2d+1)\sigma^4+\frac{(4c(1-a)+(1-ac)^2)
\sigma^2}{(1-c)^2}+\frac{c^2(1-a^2)\mu}{(1-c)^3},
\end{equation*}
from which the skewness and kurtosis indexes may be promptly expressed as 
\begin{equation*}
\alpha_3=\frac{m_3}{\sigma^3}=\frac{D}{\sigma}+\frac{d\sigma}{D}+
\frac{c(1-a)}{(1-c)^2\sigma D}
\end{equation*}
and
\begin{equation*}
\alpha_4=\frac{m_4}{\sigma^4}=3(2d+1)+\frac{4c(1-a)+(1-ac)^2}{(1-c)^2
\sigma^2}+\frac{c^2(1-a^2)}{(1-c)^3\sigma^2D}.
\end{equation*}
On the basis of these expressions, we just provide a sketch of the
flexibility of the Tempered Discrete Linnik distribution in Figures 1-3.
Indeed, these Figures display the parametric plot of $(\alpha_3,\alpha_4)$ as
$c$ and $d$ vary (for fixed $a$ and $b$) and it is apparent that the Tempered
Discrete Linnik distribution may cover an extended region in the
$(\alpha_3,\alpha_4)$-plane.

\centerline{\bf Figure 1 about here}
\centerline{\bf Figure 2 about here}
\centerline{\bf Figure 3 about here}
\centerline{\bf Figure 4 about here}

\vskip0.5cm\noindent{\bf 7. Conclusions}\vskip0.15cm\noindent We have proposed a tempered version of the Discrete Linnik
law. The new distribution is introduced on the basis of several stochastic
mixture and compound representations, which strongly justify its genesis. The
Tempered Discrete Linnik law encompasses many well-known distributions as
special cases and displays appealing probabilistic features. In addition, we
provide a manageable form for the corresponding probability function, which
may enhance the practical adoption of the Tempered Discrete Linnik law as a
statistical model for many different types of dataset.
\vskip0.5cm\noindent{\bf Acknowledgements}
The authors would like to thank Prof. Luca Pratelli for many useful advices.

\makeatletter

\renewenvironment{thebibliography}[1]
     {\noindent{\bf References}%
      \@mkboth{\MakeUppercase\bibname}{\MakeUppercase\bibname}%
      \list{}%
           {\settowidth\labelwidth{\@biblabel{#1}}%
            \leftmargin\labelwidth
            \advance\leftmargin\labelsep
            \itemindent=-\leftmargin
            \@openbib@code
            \usecounter{enumiv}%
            \let\p@enumiv\@empty
            \renewcommand\theenumiv{\@arabic\c@enumiv}}%
      \sloppy
      \clubpenalty4000
      \@clubpenalty \clubpenalty
      \widowpenalty4000%
      \sfcode`\.\@m}
     {\def\@noitemerr
       {\@latex@warning{Empty `thebibliography' environment}}%
      \endlist}

\vskip6mm

\noindent Lucio Barabesi, {\it Dipartimento di Economia Politica e Statistica},
P.zza S.Francesco 7, 53100 Siena, Italy.\\ \noindent E-mail: lucio.barabesi@unisi.it

\vskip0.3cm\noindent Carolina Becatti, {\it IMT School for advanced studies},
Piazza S.Francesco, 19, 55100 Lucca.\\ \noindent E-mail: carolina.becatti@imtlucca.it 

\vskip0.3cm
\noindent Marzia Marcheselli, {\it Dipartimento di Economia Politica e Statistica},
P.zza S.Francesco 7, 53100 Siena, Italy.\\ \noindent E-mail: marzia.marcheselli@unisi.it

\vfill\eject

\begin{figure}
\centerline{\includegraphics{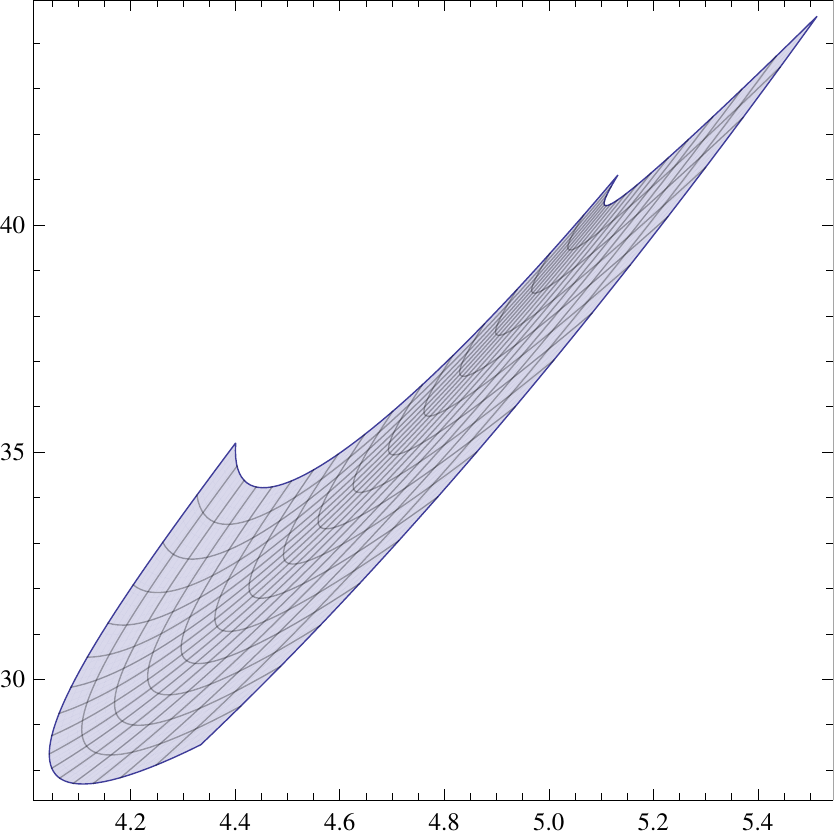}}
\caption{Parametric plot of $(\alpha_3,\alpha_4)$ for $a=1/4$,
$b=1$ and $c\in[3/10,7/10]$, $d\in[-1,3]$.}
\end{figure}

\begin{figure}
\centerline{\includegraphics{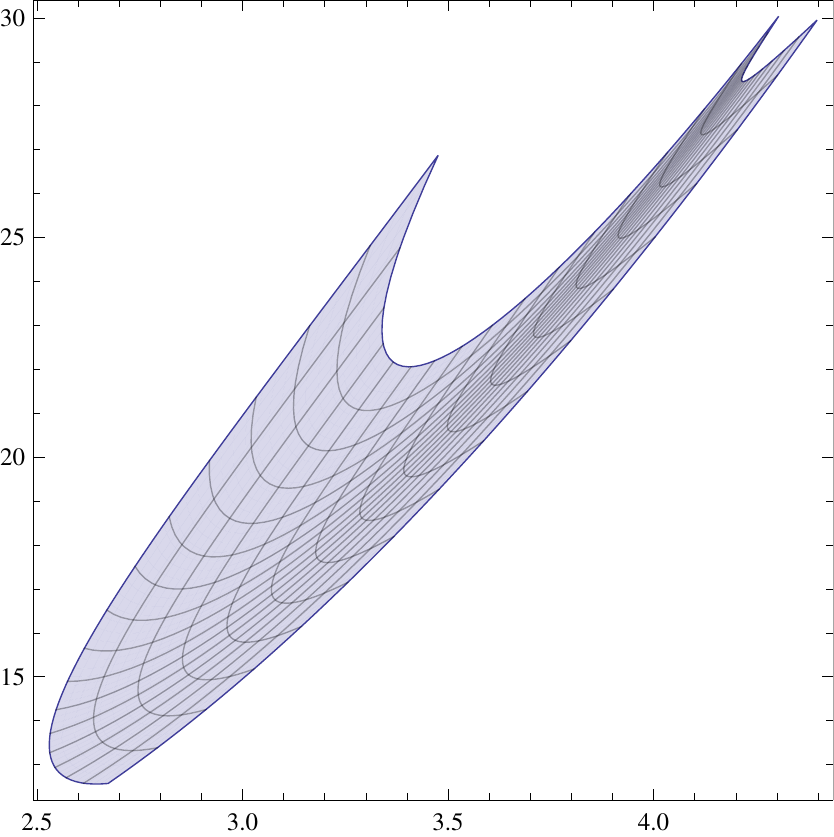}}
\caption{Parametric plot of $(\alpha_3,\alpha_4)$ for $a=1/2$,
$b=1$ and $c\in[3/10,7/10]$, $d\in[-1,3]$.}
\end{figure}

\begin{figure}
\centerline{\includegraphics{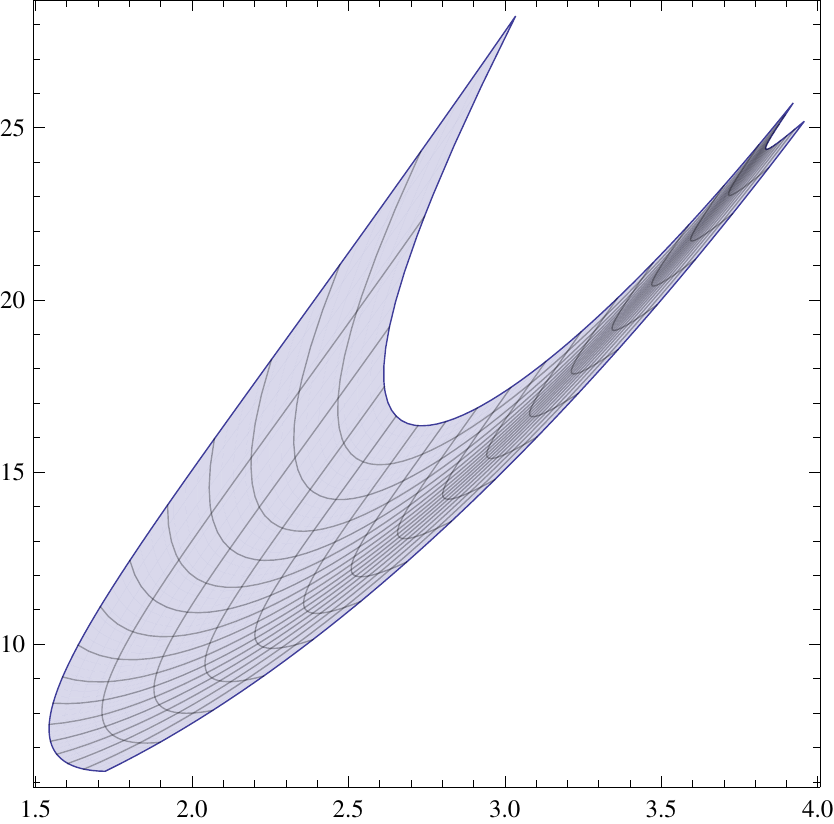}}
\caption{Parametric plot of $(\alpha_3,\alpha_4)$ for $a=3/4$,
$b=1$ and $c\in[3/10,7/10]$, $d\in[-1,3]$.}
\end{figure}

\begin{figure}
\centerline{\includegraphics{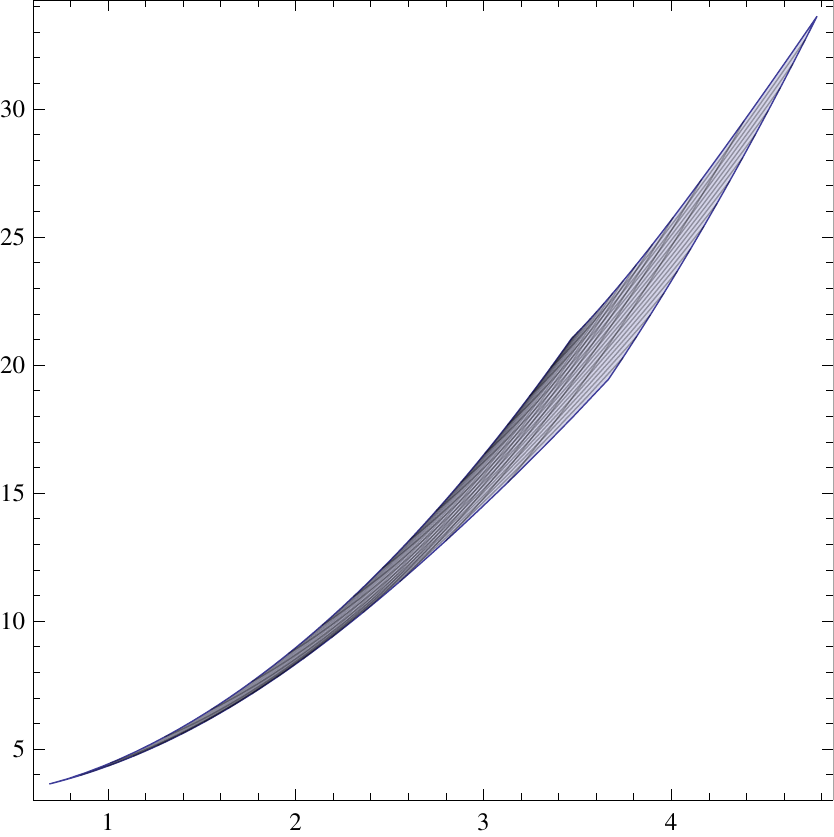}}
\caption{Parametric plot of $(\alpha_3,\alpha_4)$ for $a=-1$,
$b=1$ and $c\in[1/10,9/10]$, $d\in[0,3]$.}
\end{figure}

\end{document}